\documentclass[a4paper,10pt,reqno]{amsart}
\usepackage[a4paper,tmargin=40mm,bmargin=35mm,hmargin=30mm]{geometry}
\usepackage[T1]{fontenc}
\usepackage{lmodern}
\usepackage{amsmath}
\usepackage{amssymb}
\usepackage{amsthm}
\usepackage{stmaryrd}
\usepackage{tikz}
\usepackage[all]{xy}
\usepackage{booktabs}
\usepackage{multirow}
\usepackage{multicol}
\usepackage{here}
\usepackage{aliascnt}
\usepackage{hyperref}
\usepackage[noabbrev]{cleveref}

\newcommand{\NewTheorem}[2]{
	\newaliascnt{#1}{TheoremEnvironment}
	\newtheorem{#1}[#1]{#1}
	\aliascntresetthe{#1}
	\crefname{#1}{#1}{#2}
	\Crefname{#1}{#1}{#2}
}

\theoremstyle{definition}

\NewTheorem{Definition}{Definitions}
\NewTheorem{Remark}{Remarks}
\NewTheorem{Axiom}{Axioms}
\NewTheorem{Example}{Examples}
\NewTheorem{Observation}{Observations}
\NewTheorem{Convention}{Conventions}
\NewTheorem{Notation}{Notations}
\NewTheorem{Setting}{Settings}
\NewTheorem{Question}{Questions}
\NewTheorem{Answer}{Answers}
\NewTheorem{Conjecture}{Conjectures}
\NewTheorem{Problem}{Problems}
\NewTheorem{Solution}{Solutions}
\NewTheorem{Goal}{Goals}
\NewTheorem{Comment}{Comments}
\NewTheorem{Aim}{Aims}
\NewTheorem{Caution}{Cautions}
\NewTheorem{Exercise}{Exercises}
\NewTheorem{Examples}{Examples}
\theoremstyle{plain}
\NewTheorem{Proposition}{Propositions}
\NewTheorem{Lemma}{Lemmas}
\NewTheorem{Theorem}{Theorems}
\NewTheorem{Corollary}{Corollaries}

\crefname{enumi}{}{}
\Crefname{enumi}{}{}
\creflabelformat{enumi}{(#2#1#3)}
\crefname{enumii}{}{}
\Crefname{enumii}{}{}
\creflabelformat{enumii}{(#2#1#3)}
\crefname{enumiii}{}{}
\Crefname{enumiii}{}{}
\creflabelformat{enumiii}{(#2#1#3)}

\makeatletter
\renewcommand{\p@enumii}{}
\renewcommand{\p@enumiii}{}
\makeatother

\numberwithin{equation}{section}
\crefname{equation}{}{}
\Crefname{equation}{}{}
\creflabelformat{equation}{(#2#1#3)}

\newcommand{\SwapSymbols}[1]{
	\expandafter\let\expandafter\temporarysymbol\csname #1\endcsname
	\expandafter\let\csname #1\expandafter\endcsname\csname var#1\endcsname
	\expandafter\let\csname var#1\endcsname\temporarysymbol
}

\SwapSymbols{epsilon}
\SwapSymbols{phi}
\SwapSymbols{Gamma}
\SwapSymbols{Delta}
\SwapSymbols{Theta}
\SwapSymbols{Lambda}
\SwapSymbols{Xi}
\SwapSymbols{Pi}
\SwapSymbols{Sigma}
\SwapSymbols{Upsilon}
\SwapSymbols{Phi}
\SwapSymbols{Psi}
\SwapSymbols{Omega}

\newcommand{\bbM}{\mathbb{M}}
\newcommand{\bbN}{\mathbb{N}}

\newcommand{\cA}{\mathcal{A}}

\newcommand{\cC}{\mathcal{C}}
\newcommand{\cD}{\mathcal{D}}

\newcommand{\cF}{\mathcal{F}}

\newcommand{\cM}{\mathcal{M}}
\newcommand{\cN}{\mathcal{N}}

\newcommand{\cS}{\mathcal{S}}
\newcommand{\cT}{\mathcal{T}}

\newcommand{\cW}{\mathcal{W}}

\newcommand{\To}{\longrightarrow}

\DeclareMathOperator{\Hom}{Hom}

\DeclareMathOperator{\Ext}{Ext}
\DeclareMathOperator{\Tor}{Tor}

\DeclareMathOperator{\Mod}{Mod}

\DeclareMathOperator{\Ann}{Ann}

\DeclareMathOperator{\Coker}{Coker}

\DeclareMathOperator{\Spec}{Spec}

\DeclareMathOperator{\Min}{Min}
\DeclareMathOperator{\Ass}{Ass}
\DeclareMathOperator{\Supp}{Supp}

\title{Melkersson condition for extension of Serre subcategories}

\author{Ismael Akray, Runak H. Mustafa and Reza Sazeedeh}
\address{Department of Mathematics, Soran University, Soran, Kurdistan Region, Iraq}
\email{akray.ismael@gmail.com}
\address{Department of Mathematics, Soran University, Soran, Kurdistan Region, Iraq}
\email{rhm310h@maths.soran.edu.iq}
\address{Department of Mathematics, Urmia University, P.O.Box: 165, Urmia, Iran}
\email{rsazeedeh@ipm.ir}

\begin{document}

\begin{abstract}
Let $R$ be a commutative noetherian ring and let $\frak a$ be an ideal of $R$. In this paper, we study a certain condition, namely $C_{\frak a}$, introduced by Aghapournahr and Melkersson, on the extension of two subcategories of $R$-modules. We extend and generalize some of the main results of Yoshizawa [Y1,Y2]. As an example of extension of subcategories, we study the weakly Laskerian  modules and we find some conditions under which the local cohomology modules of a weakly Laskerian module lie in an arbitrary Serre subcategory. Eventually, we investigate the cofiniteness of the local cohomology modules of weakly Laskerian modules.   
  	\end{abstract}

\maketitle
\tableofcontents

\section{Introduction}
Throughout this paper $R$ is a commutative noetherian ring and $\frak a$ is an ideal of $R$. Let $\cS$ be a Serre subcategory of $R$-modules, let $M$ be an $R$-module and let $n$ be a non-negative integer. It is a natural question to ask:
\begin{center}
When does $H_{\frak a}^i(M)\in\cS$ for all $i\leq n$ (or $H_{\frak a}^i(M)\in\cS$ for all $i\geq n$)?
\end{center}
Aghapournahr and Melkersson gave an answer when $\cS$ satisfies the following condition:
\begin{center}
$C_{\frak a}$: If $M=\Gamma_{\frak a}(M)$ and $(0:_M\frak a)\in\cS$, then $M\in\cS$.
\end{center}

 As some of the well-known Serre subcategories are induced by extension of two subcategoris (such as the subcategory of minimax modules [Z] and the subcategory of weakly Laskerian modules [DMa, B]), Section 2 is devoted to extension of subcategories. In \cref{co2}, we give a characterization which specifies when the extension of two subcategories is closed under extensions. 

In Section 3, we want to find out when the extension of two subcategories satisfies the condition $C_{\frak a}$. Let $\cS_1$ and $\cS_2$ be subcategories of $R$-modules. The extension subcategory of $\cS_1$ and $\cS_2$ is denoted by $\cS_1\cS_2$. Assume that $\cS_1$ is Serre and $\cS1\cS_2$ is closed under extensions. In \cref{me}, we show that if $\cS_2\cS_1$ satisfies the condition $C_{\frak a}$, then so does $\cS_1\cS_2$. Assume that $\cN$ is the subcategory of finitely generated $R$-modules and $\cS$ is a Serre subcategory of $R$-modules. In \cref{hh}, we prove that if $\cS$ satisfies the condition $C_{\frak a}$ and $R/\frak a\in\cS$, then $\cN\cS$ satisfies the condition $C_{\frak a}$. The converse of this implication holds if $\frak a$ has positive height. 
 We define finitely weakly Laskrian modules and in \cref{lfk}, we prove that, $\cF w$, the subcategory of finitely weakly Laskerian modules is closed under submodules and extension. Further, we show that $\cF w$ satisfies the condition $C_{\frak a}$ for each ideal $\frak a$ of $R$. 

Fro every Serre subcategory $\cS$ of $R$-modules, Yoshizawa [Y2] studied $\bbM[\cS]$ a subset of $\Spec R$ consisting of all prime ideals $\frak p$ such that $\cS$ satisfies the condition $C_{\frak p}$. In Section 4, we generalize some of the main result of [Y2] by removing this condition that the Serre subcategories would be closed under injective envelopes. In \cref{text}, we prove that $\Supp[\cS]\cap\bbM[\cS]\subseteq\bbM[\cN\cS]$. In particular, $\Supp[\cS]_{\geq 1}\cap\bbM[\cS]=\bbM[\cN\cS]_{\geq 1}$ where $\bbM[\cS]_{\geq 1}$ consists all prime ideal $\frak p$ such that ht$\frak p\geq 1$. Let $(R,\frak m)$ be a local ring and $\frak a$ be an ideal of $R$ such that ${\rm ht} \frak a=\dim R-1$. We prove that $\bbM[\cS_{\cA.\frak a-cof}]_{\geq\dim R-1}\subseteq\Min \frak a$ and $\bbM[\cS_{\cA.\frak a-cof}]_{\geq\dim R-1}$ has at most one element if $\dim R\geq 1$ and 
$\bbM[\cS_{M.\frak a-cof}]_{\geq \dim R-1}=\emptyset$ if $\dim R\geq 2$ (cf. See \cref{ppp} and \cref{pppp}).

Section 5 is devoted to local cohomology of laskerian modules. For unexplained terminology of local cohomology, we refer to [BS]. In this section, we extend the main results of [AM] to weakly Laskerian modules (see \cref{s11} and \cref{serw}). In the rest of this section, we study the cofinitness of local cohomology of Weakly Laskerian modules. An $R$-module $M$ is said to be $\frak a$-{\it cofinite} if $\Supp M\subseteq V(\frak a)$ and
$\Ext_R^i(R/\frak a, M)$ is finite for all integers $i\geq 0$. This concept was defined for the first time by Hartshorne [H], giving a negative answer to a question of [G, Expos XIII, Conjecture 1.1].  We prove that if $M$ is a weakly Laskerian module and $\Ext_R^i(R/\frak a,M)$ is finitely generated for $i=0,1$, then $\Gamma_{\frak a}(M)$ is $\frak a$-cofinite (see \cref{pp}). Moreover, if $\dim R/\frak a=1$, then $H_{\frak a}^i(M)$ is $\frak a$-cofinite for all $i\geq 0$  (see \cref{ll}). For the case where $\dim R/\frak a=2$, we extend some main results of [NS] for the weakly Laskerian modules (see \cref{dimtwo} and \cref{dim3}).


\medskip
\section{Extension of subcategories}
Let $\cC$ be an abelian category and $\cS$ be a subcategory of $\cC$. We denote by $\cS_{\rm quot}$, the smallest subcategory of $\cC$ containing $\cS$ which is closed under quotients. This subcategory can be specified as follows:    

$$(\cS)_{\rm quot}=\{M\in\cC|\hspace{0.1cm} M \hspace{0.1cm}{\rm is\hspace{0.1cm} a\hspace{0.1cm} quotient\hspace{0.1cm}
object\hspace{0.1cm} of\hspace{0.1cm} an \hspace{0.1cm}object\hspace{0.1cm} of}\hspace{0.1cm} \cS\};$$

Let $\cT$ be another subcategory of $\cC$. We denote by $\cS\cT$, the extension subcategory of $\cS$ and $\cT$ which is:

$$\cS\cT=\{M\in\cC|\hspace{0.1cm} {\rm there \hspace{0.1cm} exists \hspace{0.1cm} an\hspace{0.1cm}exact\hspace{0.1cm}
 sequence}\hspace{0.1cm}
0\To L\To M\To N\To 0\hspace{0.1cm} {\rm with}\hspace{0.1cm} L\in\cS
\hspace{0.1cm}{\rm and}\hspace{0.1cm} N\in\cT\}.$$

For any $n\in \mathbb{N}_0$, we set $\cS^0=\{0\}$ and
$\cS^n=\cS^{n-1}\cS$. In the case where $\cS^2=\cS$, we say
that $\cS$ is {\it closed under extension}. We also define
$(\cS)_{\rm ext}=\bigcup_{n\geq 0}\cS^n$ the smallest
subcategory of $\cC$ containing $\cS$ which is closed under
extension.

A full subcategory $\cS$ of $\cC$ is called {\it Serre} if it is
closed under taking subobjects, quotients and extensions. 

\medskip

\begin{Lemma}\label{ass}
Let $\cS_1,\cS_2$ and $\cS_3$ be subcategories of an abelian category $\cC$. Then $\cS_1(\cS_2\cS_3)=(\cS_1\cS_2)\cS_3$.
\end{Lemma}
\begin{proof}
Given $M\in\cS_1(\cS_2\cS_3)$, there exists some exact sequences $0\To S_1\To M\To S\To 0$ and $0\To S_2\To S\To S_3\To 0$ such that $S_i\in\cS_i$ for $1\leq i\leq 3$. The pullback diagram of these exact sequences yields an exact sequence $0\To Y\To M\To S_3\To 0$ such that $Y\in\cS_1\cS_2$; and hence $M\in (\cS_1\cS_2)\cS_3$. The other inclusion is obtained similarly.
\end{proof}

\medskip

\begin{Proposition}\label{exte}
Let $\cS$ and $\cT$ be subcategories of an abelian category $\cC$ containing zero object. Then the following conditions hold.

${\rm (a)}$ $(\cS\cT)_{\rm sub}\subseteq (\cS)_{\rm sub}(\cT)_{\rm sub}.$

${\rm (b)}$ $(\cS\cT)_{\rm quot}\subseteq (\cS)_{\rm quot}(\cT)_{\rm quot}.$

${\rm (c)}$ The following conditions are equivalent.

${\rm(i)}$ $(\cS\cT)_{\rm ext}\subseteq (\cS)_{\rm ext}(\cT)_{\rm ext}$.

 ${\rm(ii)}$ $(\cT\cS)_{\rm ext}\subseteq (\cS)_{\rm ext}(\cT)_{\rm ext}$.

If $\cT$ is closed under extensions, then ${\rm (i), (ii)}$ are equivalent to:

 ${\rm(iii)}$ $\cT\cS\subseteq (\cS)_{\rm ext}\cT$.
\end{Proposition}
\begin{proof}
(a),(b) follows from [K, Proposition 2.4]. (c)  In order to prove, (i)$\Rightarrow$(ii), in view of \cref{ass}, for every positive integer $n$, we have $(\cT\cS)^n=\cT(\cS\cT)^{n-1}\cS\subseteq (\cS\cT)_{\rm ext}$ as $\cT,(\cS\cT)^{n-1}$ and $\cS$ are subclasses of $(\cS\cT)_{\rm ext}$. (ii)$\Rightarrow$(i). For every positive integer $n$, we show that $(\cS\cT)^n\subseteq (\cS)_{\rm ext}(\cT)_{\rm ext}$. The assertion for $n=1$ is clear. For $n\geq 2$, it follows from \cref{ass} that $(\cS\cT)^n=\cS(\cT\cS)^{n-1}\cT$.
Then for every $M\in(\cS\cT)^n$, there exist an exact sequence $0\To N\To M\To K\To 0$  such that $N\in\cS(\cT\cS)^{n-1}$ and $K\in\cT$. Furthermore, there exist an exact sequence $0\To S\To N\To X\To 0$ such that $S\in\cS$ and $X\in(\cT\cS)^{n-1}$. By the assumption, there exist positive numbers $r,t$ such that $X\in\cS^r\cT^t$ so that $N\in\cS^{r+1}\cT^t$. Therefore, $M\in\cS^{r+1}\cT^{t+1}\subseteq (\cS)_{\rm ext}(\cT)_{\rm ext}$. The implication (ii)$\Rightarrow$(iii) is clear. In order to prove (iii)$\Rightarrow$(ii), we show by induction on $n$ that $(\cT\cS)^n\subseteq (\cS)_{\rm ext}\cT$. The case $n=1$ is the assumption. Assume that $n>1$ and the result has been proved for all values smaller that $n$. In view of \cref{ass}, we have $(\cT\cS)^n=(\cT\cS)^{n-1}\cT\cS\subseteq(\cS)_{\rm ext}\cT\cT\cS=(\cS)_{\rm ext}\cT\cS\subseteq (\cS)_{\rm ext}(\cS)_{\rm ext}\cT=(\cS)_{\rm ext}\cT$. 
\end{proof}
  
\medskip
\begin{Corollary}\label{co2}
Let $\cS$ and $\cT$ be subcategories of an abelian category which are closed under extensions. Then $\cS\cT$ is closed under extension if and only if $\cT\cS\subseteq \cS\cT$. 
\end{Corollary}
\begin{proof}
Straightforward by \cref{exte}.
\end{proof}

We denote by $\cN$, the Serre subcategory of finitely generated $R$-modules. We have the following corollary.
\medskip
\begin{Corollary}\label{ccc}
If $\cS$ is a subcategory of $R$-modules, then $\cS\cN\subseteq \cN(\cS)_{\rm quot}$. In particular if $\cS$ is closed under quotients and extensions, then $\cN\cS$ is closed under quotients and extensions.
\end{Corollary}
\begin{proof}
Given $M\in\cS\cN$, there exists an exact sequence of $R$-modules  $0\To S\To M\To K\To 0$ such that $S\in\cS$ and $K$ is finitely generated. Clearly, there exists a finitely generated submodule $N$ of $M$ such that $M=S+N$ and hence $M$ is a quotient of $N\oplus S\in \cN\cS$ so that $M\in(\cS\cN)_{\rm quot}\subseteq \cN(\cS)_{\rm quot}$ by \cref{exte}(b). The second assertion follows by the first assertion and \cref{co2}. 
\end{proof}

\medskip

\begin{Example}
Let $\cD=\{M\in R$-Mod| Att$_RM$ is a finite set$\}$. By the basic properties of attached primes, $\cD$ is closed under quotients and extensions; and hence it follows from \cref{ccc} that $\cN\cD$ is closed under quotients and extensions.
\end{Example}
\section{Extension of subcategories and Mekersson condition}

Throughout this section $\frak a$ is an ideal of $R$. Given a subcategory $\cS$ of modules and an $R$-module $M$, the subcategory $\cS$ is said to satisfy {\it the condition $C_{\frak a}$
on} $M$ if the following implication holds:
\begin{center}
If $\Gamma_{\frak a}(M)=M$ and $(0:_M{\frak a})$ is in $\cS$, then
$M$ is in $\cS$.
\end{center}
We denote by $\cS_{\frak a}$, the largest subcategory of modules such
that $\cS$ satisfies the condition $C_{\frak a}$ on $\cS_{\frak a}$.
Clearly, $\cS\subseteq\cS_{\frak a}$. The class $\cS$ is said to satisfy {\it the condition} $C_{\frak a}$
whenever $\cS_{\frak a}= R$-Mod.

\medskip

\begin{Proposition}\label{sed}
Let $\cS$ be a subcategory of modules which is closed under, submodules, extensions and direct unions. Then it satisfies $C_{\frak a}$ condition for every ideal $\frak a$ of $R$.
\end{Proposition}
\begin{proof}
Let $\frak a$ be an ideal of $R$ and $M$ be an $R$-module such that $M=\Gamma_{\frak a}(M)$ and $(0:_M\frak a)\in\cS$. We prove that $(0:_M\frak a^n)\in\cS$ for all $n\geq 1$. It suffices to show that for $n=2$. There exist the exact sequences $0\To \frak a/\frak a^2\To R/\frak a^2\To R/\frak a\To $ and $0\To X \To (R/\frak a)^t\To \frak a/\frak a^2\To 0$ for some positive integer $t$. Applying the functor $\Hom_R(-,M)$ and using the fact that $\cS$ is closed under submodules and extensions, we deduce that $(0:_M\frak a^2)\in\cS$. Since $\cS$ is closed under direct unions, $M=\Gamma_{\frak a}(M)\in\cS$.
\end{proof}

\begin{Example}
We notice that the converse of \cref{sed} may not be hold. To be more precise, if we consider $\cD=\{M\in R$-Mod$| \hspace{0.1cm}\Ass _RM$ is a finite set$\}$, then it is clear that $\cD$ satisfies the condition $C_{\frak a}$ for each ideal $\frak a$ of $R$ while $\cD$ is not closed under direct unions. 
\end{Example}
\medskip
The following proposition generalizes [Y1, Theorem 4.3].

\begin{Proposition}\label{me}
Let $\frak a$ be an ideal of $R$, let $\cS_1$ be a Serre subcategory of $R$-modules and let $\cS_2$ be any subcategory of modules containing the zero module. If $\cS_1\cS_2$ is closed under extensions and $(\cS_2\cS_1)_{\frak a}$ is closed under quotients, then $(\cS_2\cS_1)_{\frak a} \subseteq(\cS_1\cS_2)_{\frak a}$. In particular, if $\cS_2\cS_1$ satisfies $C_{\frak a}$ condition, then so is $\cS_1\cS_2$.
\end{Proposition}
\begin{proof}
Given $M\in (\cS_2\cS_1)_{\frak a}$ with $M=\Gamma_{\frak a}(M)$ and $(0:_M\frak a)\in\cS_1\cS_2$, there exist an exact sequence $$0\To L_1\To (0:_M\frak a)\To L_2\To 0$$ such that $L_i\in\cS_i$ for $i=1,2$. We now have a pushout diagram 
$$\xymatrix{&&0\ar[d]&0\ar[d]\\
0\ar[r]&L_1\ar[r]\ar@{=}[d]& (0:_M\frak a)\ar[d]\ar[r]& \ar[d]L_2\ar[r]\ar[d]\ar[r]& 0\\
0\ar[r]&L_1\ar[r]& M\ar[r]\ar[d]& D\ar[d]\ar[r]& 0\\
 &&\overline{M}\ar[d]\ar@{=}[r]& \overline{M}\ar[d]\\
&&0&0.}$$ 
Application of the functor $\Hom_R(R/\frak a,-)$ to the middle row, induces the following exact sequence $$0\To L_1\To (0:_M\frak a)\To (0:_D\frak a)\To \Ext_R^1(R/\frak a,L_1).$$ It is clear that $\Ext^1(R/\frak a,L_1)\in\cS_1$ so that $(0:_D\frak a)\in\cS_2\cS_1$. Now, since by the assumption $D\in(\cS_2\cS_1)_{\frak a}$, we get that $D\in\cS_2\cS_1$. Thus, according to \cref{exte} and the fact that $\cS_1\cS_2$ is closed under extension, we have $M\in\cS_1(\cS_2\cS_1)=(\cS_1\cS_2)\cS_1\subseteq \cS_1\cS_2$.
\end{proof}

\medskip
\begin{Corollary}\label{co4}
Let  $\frak a$ be an ideal of $R$, let $\cS_1$ be a Serre subcategory of $R$-modules and let $\cS_2$ is closed under extensions such that $\cS_2\cS_1\subseteq \cS_1\cS_2$. If $\cS_2\cS_1$ satisfies the condition $C_{\frak a}$, then so does $\cS_1\cS_2$. 
\end{Corollary}
\begin{proof}
It follows from \cref{co2} that $\cS_1\cS_2$ is closed under extension. Now the result is obtained by \cref{me}. 
\end{proof}

\medskip
\begin{Corollary}
Let $\frak a$ be an ideal of $R$, let $\cS$ be a subcategory of $R$-modules closed under quotients and extensions. If $\cS\cN$ satisfies the condition $C_{\frak a}$, then so does $\cN\cS$. 
\end{Corollary}
\begin{proof}
It follows from \cref{ccc} that $\cS\cN\subseteq \cN\cS$. Now the result follows from \cref{co4}.
\end{proof}

\medskip
For any Serre subcategory $\cS$ of modules we set $\Supp\cS=\{\frak p\in\Spec R|\hspace{0.1cm} R/\frak p\in\cS\}$
\medskip
\begin{Proposition}\label{teh1}
Let $\frak a$ be an ideal of positive height and let $\cS$ be a Serre subcategory of modules. If $(\cN\cS)_{\frak a}\subseteq \cS_{\frak a}$. In particular, if $\cN\cS$ satisfies the condition $C_{\frak a}$, then so does $\cS$.
\end{Proposition}
\begin{proof}
Given $M\in\cN\cS$ with $M=\Gamma_{\frak a}(M)$ and $(0:_M\frak a)\in\cS$, there exists an exact sequence $0\To N\To M\To S\To 0$ of $R$-modules such that $N$ is finitely generated and $S\in\cS$.  It is clear that $\Supp N=\Supp (0:_N\frak a)\subseteq \Supp(0:_M\frak a)\subseteq \Supp\cS$. Considering a finite filtration $0=N_0\subseteq N_1\subseteq \dots \subseteq N_n=N$ for $N$ with prime ideals $\frak p_i$ such that $N_i/N_{i-1}\cong R/\frak p_i\in\cS$, $1\leq i\leq n$, we deduce that $N\in\cS$; and hence $M\in\cS$. It now follows from [SR, Proposition 2.3] that $(\cN\cS)_{\frak a}\subseteq \cS_{\frak a}$. The second assertion is clear by the first.
\end{proof}

\medskip

\begin{Lemma}\label{iii}
Let $\frak a$ be an ideal of $R$ of positive height and let $\cS$ be a Serre subcategory of modules such that $R/\frak a\notin\cS$. Then there exists $\frak p\in\Min \frak a$ such that $H_{\frak a}^{{\rm ht}\frak p}(R)\notin\cN\cS$. 
\end{Lemma}
\begin{proof}
It is clear that $R/\frak a\notin\cS$ if and only if $R/\sqrt{\frak a}\notin\cS$ and so we may assume that $\frak a=\sqrt{\frak a}$. Assume that $H_{\frak a}^{\rm ht\frak p}(R)\in\cN\cS$ for all $\frak p\in\Min \frak a$. Since $R/\frak a\notin\cS$, there exists $\frak p\in\Min\frak a$ such that $R/\frak p\in\cS$. Assume that ${\rm ht}\frak p=t$ and we claim that $H_{\frak a}^t(R)\notin\cN\cS$. Otherwise there exists an exact sequence $0\To F\To H_{\frak a}^t(R)\To S\To 0$ such that $F$ is finitely generated and $S\in\cS$. Localizing in $\frak p$, since $\cS$ is Serre and $R/\frak p\notin\cS$, we have $S_{\frak p}=0$. This implies that $F_{\frak p}\cong H_{\frak pR_{\frak p}}^t(R_{\frak p})$ is finitely generated which contradicts the non-vanishing theorem in local cohomology (see [BS, 6.1.7 Exercise]).
\end{proof}

\medskip
\begin{Theorem}\label{hh}
Let $\frak a$ be an ideal of $R$ and let $\cS$ be a Serre subcategory of $R$-modules. If $\cS$ satisfies the condition $C_{\frak a}$ and $R/\frak a\in\cS$, then $\cN\cS$ satisfies the condition $C_{\frak a}$. The converse of this implication holds if $\frak a$ has positive height. 
\end{Theorem}
\begin{proof}
Given an $R$-module $M$ with $M=\Gamma_{\frak a}(M)$ and $(0:_M\frak a)\in\cN\cS$, there exists an exact sequence $0\To N\To (0:_M\frak a)\To S\To 0$ such that $N$ is finitely generated and $S\in\cS$.
We observe that $\Supp N\subseteq V(\frak a)$ and since $N$ is noetherian, there exists a finite filtration $0=N_0\subseteq N_1\subseteq \dots \subseteq N_n=N$ such that $N_i/N_{i-1}\cong R/\frak p_i\in\cS$ for $1\leq i\leq n$. Thus $N\in\cS$ and so $(0:_M\frak a)\in\cS$. The assumption implies that $M\in\cS$; and hence $M\in\cN\cS$. To prove, the converse, it follows from \cref{teh1} that $\cS$ satisfies the condition $C_{\frak a}$. We observe that $\Ext_R^i(R/\frak a,R)\in\cN\cS$ for all $i\geq 0$; and hence it follows from [AM, Theorem 2.9] that $H_{\frak a}^i(R)\in\cN\cS$ for all $i\geq 0$. Now \cref{iii} implies that $R/\frak a\in\cS$. 
\end{proof}

We denote by $\cF$, a Serre subcategory consisting of all $R$-modules of finite support. We have the following lemma.  
\medskip
\begin{Lemma}\label{llf}
The Serre subcategory $\cF$ satisfies the condition $C_{\frak a}$ for all ideals $\frak a$ of $R$.
\end{Lemma}
\begin{proof}
Given an ideal $\frak a$ of $R$ and an $R$-module $M$ with $M=\Gamma_{\frak a}(M)$ and $(0:_M\frak a)\in F$, we have  $\Supp M=\Supp(0:_M\frak a)$; and hence $M\in\cF$.
\end{proof}

\medskip

An $R$-module $M$ is said to be {\it weakly Laskerian} if the set of associated prime ideals of any quotient module of $M$ is finite. We denote, by $\cW$, the subcategory of weakly Laskerian modules. It follows from [DMa, Lemma 2.2] that $\cW$ is a Serre subcategory of $\Mod R$. For the subcategory of weakly Laskerian modules, Bahmanpour [B] gave the following result that we give a slightly alternative short proof.

\begin{Proposition} ([B, Theorem 3.3])\label{bah} 
Let $R$ be a noetherain ring. Then $\cW=\cN\cF$.
\end{Proposition}
\begin{proof}
The assertion $\cN\cF\subseteq \cW$ is obtained as [B, Theorem 3.3]. Assume that $M\in\cW$ and suppose on the contrary that $M\notin\cN\cF$. The for any finitely generated submodule $N$ of $M$, the set $\Supp M/N$ is infinite and since $R$ is noetherian, there are infinitely many prime ideals in $\Supp M/N$ which are pairwise indecomposable under inclusions. By a similar argument in the proof of [B, Theorem 2.2], there exists a submodule $L$ of $M$ such that $\Ass M/L$ is an infinite set which is a contradiction. 
\end{proof}

\medskip

\begin{Corollary}\label{coo}
If $\frak a$ is an ideal of $R$ such that $\Supp R/\frak a$ is finite, then $\cW$ satisfies the condition $C_{\frak a}$ . In particular, if $\dim R/\frak a\leq 1$, then $\cW$ satisfies the condition $C_{\frak a}$.
\end{Corollary} 
\begin{proof}
By \cref{bah}, we have $\cW=\cN\cF$; and hence the result follows from \cref{hh}.
\end{proof}

An $R$-module $M$ is said to be {\it finitely weakly Laskerian} if $\Ass M/K$ is a finite set for any finite submodule $K$ of $M$. We denote by $\cF w$, the full subcategory of finitely weakly Laskerian modules. 

\medskip

\begin{Lemma}\label{lfk}
The following condition hold.

${\rm (i)}$ The quotient module $M/K$ is in $\cF w$ for every module $M\in\cF w$ and every finitely generated submodule $K$ of $M$.

 ${\rm (ii)}$ The subcategory $\cF w$ is closed under submodules and extensions. 
\end{Lemma}
\begin{proof}
(i) is straightforward. (ii) It is clear that $\cF w$ is closed under submodules. Assume that $0\To K\To M\To L\To 0$ is an exact sequence such that $K,L\in\cF w$ and consider a finitely generated submodule $X$ of $M$. Considering $K$ as a submodule of $M$, if $X\cap K=0$, we have an exact sequence $0\To K\To M/X\To\frac{M}{X+K}\To 0$ of modules. We observe that $\frac{X+K}{K}$ is a finitely generated submodule of $L$ so that $\Ass \frac{M}{X+K}$ is finite; and hence $\Ass M/X$ is finite. If $Y=X\cap K\neq 0$, putting $L_1=X/Y$, there exists the following commutative diagram 
$$\xymatrix{
0\ar[r]&Y\ar[r]\ar[d]^{\alpha}& X\ar[d]^{\beta}\ar[r]& L_1\ar[r]\ar[d]^{\gamma}\ar[r]& 0\\
0\ar[r]&K\ar[r]& M\ar[r]& L\ar[r]& 0}$$ which yields an exact sequence $0\To \Coker\alpha\To\Coker\beta\To \Coker\gamma\To 0 $. It follows from (i) that $\Coker\alpha, \Coker\gamma\in\cF w$ so that $\Ass\Coker\alpha$ and $\Ass\Coker\gamma$ are finite sets. This forces that $\Ass\Coker\beta$ is a finite set.
\end{proof}

\medskip
\begin{Corollary}
If $M\in(\cF w)_{\rm quot}$, then  $M$ is a direct limit of finitely weakly laskerian modules. 
\end{Corollary}
\begin{proof}
There exists a module $L\in\cF w$ and a submodule $K$ of $L$ such that
$M=L/K$. It is known that $K$ is the direct union of its finitely generated submodules, which is $K=\underset{\rightarrow}{\rm lim} K_i$. For each $i$, there is an exact sequence $0\To K_i\To L\To L/K_i\To 0$; and hence applying $\underset{\rightarrow}{\rm lim}(-)$, we deduce that $M=\underset{\rightarrow}{\rm lim}\frac{L}{k_i}$. By virtue of \cref{lfk} every $L/K_i$ is finitely weakly Laskerian. 
\end{proof}

\medskip

\begin{Proposition}\label{ls}
Let $\cF w$ be the subcategory of finitely weakly Laskerian modules. Then $\cF w$ satisfies $C_{\frak a}$ condition for each ideal $\frak a$ of $R$.
\end{Proposition}
\begin{proof}
Let $\frak a$ be an ideal of $R$ and $M$ be an $R$-module such that $M=\Gamma_{\frak a}(M)$ and $(0:_M\frak a)\in\cF w$.  Given a finitely generated submodule $K$ of $M$, there is an exact sequence of modules $0\To K\To M\To M/K\To 0$. Applying the functor $\Hom_R(R/\frak a,-)$ we have the following exact sequence of modules $0\To (0:_K\frak a)\To (0:_M\frak a)\To (0:_{M/K}\frak a)\To \Ext^1_R(R/\frak a,K)$. It follows from \cref{lfk} that $\Coker((0:_K\frak a)\To (0:_M\frak a))\in\cF w$ and $\Ext^1_R(R/\frak a,K)\in\cF w$ so that $\Coker((0:_M\frak a)\To (0:_{M/K}\frak a))\in\cF w$. This forces $\Ass((0:_{M/K}\frak a))$ is a finite set. On the other hand, it is clear that $(0:_{M/K}\frak a)$ is an essential submodule of $\Gamma_{\frak a}(M/K)=M/K$, and hence $\Ass((0:_{M/K}\frak a))=\Ass M/K$ is a finite set. 
\end{proof}
\medskip


\section{Melkersson subsets of extension of Serre subcategories}

In this section we study the Melkersson subset induced by Yoshizawa [Y2] and we generalize and extend some of the main results of [Y2] which has been proved for a prime ideal and for those Serre subcategories which are closed under injective envelopes. Let $\cS$ be a subcategory of modules. We first recall from [Y2] the following definition.

\begin{Definition}
The {\it Melkersson subset} $\bbM[\cS]$ of $\Spec R$ consists of all prime ideals $\frak p$ such that $\cS$ satisfies $C_{\frak p}$ condition. For any non-negative integer $i$, the subset $\bbM[\cS]_{\geq i}$ consists of all $\frak p\in\bbM[\cS]$ such that ht$\frak p\geq i$.  We also recall that $\cS$ is a {\it Melkersson subcategory} provided that $\bbM[\cS]=\Spec R$.
\end{Definition}
 We notice that a Melkersson subcategory of $R$-modules is not a Serre subcategory in general. To be more precise, $\cF w$ is a Melkersson subcatgory of $R$-modules while it is not Serre. Yoshizawa [Y2, Theorem 3.4] proved the following result  for those Serre subcategoris of $R$-modules which are closed under injective envelopes. Here we remove this condition and give a more general version.
\medskip

\begin{Theorem}\label{text}
Let $\cS$ be a Serre subcategory of modules. Then $\Supp[\cS]\cap\bbM[\cS]\subseteq\bbM[\cN\cS]$. In particular,
$\Supp[\cS]_{\geq 1}\cap\bbM[\cS]=\bbM[\cN\cS]_{\geq 1}$. 
\end{Theorem}
\begin{proof}
Given $\frak p\in\Supp[\cS]\cap\bbM[\cS]$ and an $R$-module $M$ with $\Gamma_{\frak p}(M)=M$ and $(0:_M\frak p)\in\cN\cS$, there exists an exact sequence $0\To N\To (0:_M\frak p)\To S$ such that $N$ is finitely generated and $S\in\cS$. Since $\cS$ is Serre, $R/\frak p\in\cS$ and $N$ is noetherian, there exist prime ideals $\frak p_i\in\Supp N$, $1\leq i\leq n$ and a filtration $0=N_0\subseteq N_1\subseteq N_1\dots \subseteq N_n=N$ such that $N_i/N_{i-1}\cong R/\frak p_i$. The fact that $\Supp N\subseteq V(\frak p)$ forces $N\in\cS$; and hence $(0:_M\frak p)\in\cS$. Since $\frak p\in\bbM[\cS]$, we deduce that $M\in\cS$ so that $M\in\cN\cS$. Conversely, according to [Y2, Proposition 3.2], we have $\bbM[\cN\cS]_{\geq 1}\subseteq \Supp[\cS]_{\geq 1}$ and so it suffices to prove that $\bbM[\cN\cS]_{\geq 1}\subseteq \bbM[\cS]$. Assume that $\frak p\in\bbM[\cN\cS]_{\geq 1}$ and $M$ is an $R$-module such that $M=\Gamma_{\frak p}(M)$ and $(0:_M\frak p)\in\cS$. Then $M\in\cN\cS$ so that there is an exact sequence $0\To N\To M\To S\To 0$ such that $N$ is finitely generated and $S\in\cS$. The fact that $\bbM[\cN\cS]\Supp\cS$ implies that $R/\frak p\in\cS$ and Since $\Supp N\subseteq V(\frak p)$, a similar argument mentioned above, deduce that $N\in\cS$ so that $M\in\cS$. 
\end{proof}

\medskip
We recall from [Z] that an $R$-module $M$ is a {\it minimax module} if there is a finitely generated submodule $N$ of $M$ such that $M/N$ is artinian. The subcategory of minimax modules is denoted by $\cM$, and so clearly $\cM=\cN\cA$ where $\cA$ is the subcategory of Artinian modules.

\begin{Corollary}
Let $\cM$ be the subcategory of minimax modules of $\Mod R$. Then $\cM$ satisfies $C_{\frak m}$ condition for each maximal ideal $\frak m$ of $R$.
\end{Corollary}
\begin{proof}
We observe that $\cM=\cN\cA$ where $\cA$ is the subcategory of artinian modules. It is known by a result of Melkersson that $\bbM[\cA]=\Spec R$. If $\dim R=0$, then $\cM=\cA$ and there is nothing to prove. If $\dim R>0$, the result follows by \cref{text}.
\end{proof}

The following result establishes a relation between the injective indecomposable modules of $\cN\cS$ and $\cS$ where $\cS$ is a  Serre subcategory of modules.

\medskip
\begin{Proposition}
Let $\cS$ be a Serre subcategory of modules and $\frak p$ be an prime ideal of $R$ such that $E(R/\frak p)\in\cN\cS$. Then $E(R/\frak p)\in\cS$ if either ${\rm ht} \frak p>0$ or ${\rm ht} \frak p=0$ and $R/\frak p\in\cS$.
\end{Proposition}
\begin{proof}
There exists an exact sequence of modules $0\To N\To E(R/\frak p)\To S\To 0$ such that $N$ is finitely generated and $S\in\cS$. For the first condition, if $\frak p\notin\Supp S$, we have $N_{\frak p}\cong E(R/\frak p)$ and hence $E(R/\frak p)$ is a finitely generated $R_{\frak p}$-module so that $R_{\frak p}$ is Cohen-Macaulay and hence $\dim R_{\frak p}=0$ which is a contradiction. Then $R/\frak p\in\cS$. Given every $\frak q\in V(\frak p)\setminus\{\frak p\}$, since $\frak p\in\Supp S$, we have $R/\frak q\in\cS$.  Thus $R/\frak q\in\cS$ for any $\frak q\in\Supp N$ as $\Ass N=\{\frak p\}$. Considering a finite filtration of 
$0=N_n\subset N_{n-1}\subset \dots\subset N_1\subset N_0=N$ of submodules of $N$ such that $N_i/N_{i+1}\cong R/\frak p_i$ for $0\leq i\leq n-1$, we deduce that $N\in\cS$. Now, the above exact sequence implies that $E(R/\frak p)\in\cS$. If the second condition holds, the previous argument conclude that $E(R/\frak p)\in\cS$.
\end{proof}

\medskip
We notice that in the above result, if $\frak p$ is a minimal ideal of $R$ such that $E(R/\frak p)\in\cN\cS$, then $E(R/\frak p)$ may not be in $\cS$ by the following example. 
\begin{Examples}
Let $(R,\frak m)$ be a $1\leq d$-dimensional Gorenstein local ring and for every non-negative integer $n$, let $\cD_n$ denote the subcategory of all modules of dimension $\leq n$. It is clear that $\cD_n$ is Serre and if $0\To R\To \bigoplus_{{\rm ht}\frak p=0} E(R/\frak p)\To X\To 0$ is an exact sequence of modules, then $X$ is a submodule of $\bigoplus_{{\rm ht}\frak p=1} E(R/\frak p)$. Then $X\in\cD_{d-1}$ and so $E(R/\frak p)\in\cN\cD_{d-1}$ for every $\frak p\in\Min R$ while $E(R/\frak p)\notin\cD_{d-1}$ for every $\frak p\in\Min R$. Specially, if $d=1$, then $E(R/\frak p)\in\cN\cA$, but $E(R/\frak p)\notin\cA$ for every $\frak p\in\Min R$ where $\cA$ is the subcategory of artinian modules. Furthermore, if $d=2$, then $E(R/\frak p)\in\cW=\cN\cF$ while $E(R/\frak p)\notin\cF$ for every $\frak p\in\Min R$.
\end{Examples} 
\medskip

 \medskip
Melkersson [M1, Theorem 1.6 and Corollary 1,7] showed that if $R$ is a local ring and $\frak a$ is an ideal of $R$, then an artinian module $M$ is $\frak a$-cofinite if and only if $(0:_M\frak a)$ has finite length and the subcategory of artinian $\frak a$-cofinite modules, denoted by $\cS_{\cA.\frak a-cof}$ is Serre. 
The following result was proved by Yoshizawa for a prime ideal (see [Y2, Proposition 6.2]) that we prove it for an arbitrary ideal $\frak a$ of $R$.

\begin{Proposition}\label{ppp}
Let $(R,\frak m)$ be a local ring of $\dim R\geq 1$ and $\frak a$ be an ideal of $R$ such that ${\rm ht} \frak a=\dim R-1$. Then the following condition holds.

${\rm (i)}$  $\cS_{\cA.\frak a-cof}=\bigcap_{\frak p\in\Min \frak a}\cS_{\cA.\frak p-cof}$.

${\rm (ii)}$ $\bbM[\cS_{\cA.\frak a-cof}]_{\geq\dim R-1}\subseteq\Min \frak a.$ Moreover, the set $\bbM[\cS_{\cA.\frak a-cof}]_{\geq\dim R-1}$ contains at most one element. 
\end{Proposition}
\begin{proof}
In order to prove (i), given an $R$-module $M\in\cS_{\cA.\frak a-cof}$, it is clear that $(0:_M\frak p)$ has finite length so that $M\in\cS_{\cA.\frak p-cof}$ for every $\frak p\in\Min\frak a$ by [M1, Theorem 1.6]. Conversely, given $M\in\cS_{\cA.\frak p-cof}$ for every $\frak p\in\Min \frak a$, it is clear that $\Gamma_{\frak a}(M)=M$ and so it follows from [DM, Corollary 1] that $M$ is $\frak a$-cofinite. To do (ii), by a similar proof of [Y2, Proposition 6.2] we show that $E(R/\frak m)\notin \cS_{\cA.\frak a-cof}$. Otherwise, $(0:_{E(R/\frak m)}\frak a)\cong E_{R/\frak a}(R/\frak m)$ is finitely generated $R/\frak a$-module so that $\dim R/\frak a=0$. But this implies that ht $\frak a=\dim R$ which is a contradiction. We now prove that $\cS_{\cA.\frak a-cof}$ satisfies the condition $C_{\frak a}$. Assume that $M$  is an $R$-module such that $\Gamma_{\frak a}(M)=M$ and $(0:_M\frak a)\in\cS_{\cA.\frak a-cof}$. Then $(0:_M\frak a)$ has finite length; and hence [M1, Theorem 1.6] implies that $M\in\cS_{\cA.\frak a-cof}$. We show that $\frak m\notin \bbM[\cS_{\cA.\frak a-cof}]$; otherwise since $(0:_{E(R/\frak m)}\frak m)$ has finite length, $(0:_{E(R/\frak m)}\frak m)\in\cS_{\cA.\frak a-cof}$; and hence $E(R/\frak m)\in\cS_{\cA.\frak a-cof}$ which is a contradiction. If $\frak q\in\bbM[\cS_{\cA.\frak a-cof}]_{\dim R-1}$, then we have $\frak a\subseteq \frak q$; otherwise it follows from [SR, Corollary 2.10] that $\cS_{\cA.\frak a-cof}$ satisfies the condition $C_{\frak a+\frak q}$ and since $\frak a+\frak q$ is $\frak m$-primary, using [SR, Proposition 2.4], we have $\frak m\in \bbM[\cS_{\cA.\frak a-cof}]$ which is a contradiction. Clearly, $\frak q\in\Min \frak a$. For the second assertion, if $\frak p$ and $\frak q$ are distinct prime ideals in $\bbM[\cS_{\cA.\frak a-cof}]_{\geq\dim R-1}$, by a similar argument $\frak p+\frak q$ is $\frak m$-primary so that $\frak m \in\bbM[\cS_{\cA.\frak a-cof}]_{\geq\dim R-1}$ which is a contradiction.
\end{proof}

 \medskip

For an ideal $\frak a$ of $R$, Melkersson [M2, Corollary 4.4] proved that the subcategory of $\frak a$-cofinite minimax modules is Serre and we denote this subcategory by $\cS_{\cM.\frak a-cof}$. The following proposition generalizes [Y2, Proposition 6.7] for an arbitrary ideal of $R$. 

\begin{Proposition}\label{pppp}
Let $(R,\frak m)$ be a local ring such that $\dim R\geq 2$ and let $\frak a$ be an ideal of $R$ such that ${\rm ht \frak a}=\dim R-1$. Then $\bbM[\cS_{M.\frak a-cof}]_{\geq \dim R-1}=\emptyset$. Furthermore, $\cS_{M.\frak a-cof}$ does not satisfy the condition $C_{\frak a}$.   
\end{Proposition}
\begin{proof}
As $\cS_{\cM.\frak a-cof}$, contains all modules of finite length, a similar argument in the proof of \cref{ppp} implies that $\frak m\notin\bbM[\cS_{\cM,\frak a-cof}]$. We now claim that $\frak p\notin\bbM[\cS_{\cM,\frak a-cof}]$ for every $\frak p\in\Min\frak a$. To do this, if $\frak p\in\bbM[\cS_{\cM.\frak a-cof}]$ for some $\frak p\in\Min \frak a$, then $\Ext_R^i(R/\frak p,R)\in\cS_{\cM.\frak a-cof}$ for each $i\geq 0$ as $V(\frak p)\subseteq V(\frak a)$. It follows from [AM, Theorem 2.9] that $H_{\frak p}^{\rm ht \frak p}(R)\in\cS_{\cM.\frak a-cof}$; and hence $H_{\frak p}^{\rm ht \frak p}(R)\in\cN\cA$. It now follows from [Y2, Lemma 3.1] that $\frak p=\frak m$ which is a contradiction. Assume that $\frak q$ be a prime ideal such that ht$\frak q=\dim R-1$. If $\frak a\subseteq \frak q$, then $\frak q\in\Min \frak a$; and hence $\frak q\notin \bbM[\cS_{\cM.\frak a-cof}]_{\geq \dim R-1}$. If $\frak a\nsubseteq \frak q$, then $\frak a+\frak p$ is $\frak m$-prmiary. Assume that $\frak q\in \bbM[\cS_{\cM.\frak a-cof}]_{\geq \dim R-1}$. We observe that $\Hom_R(R/\frak a,(0:_{E(R/\frak m)}\frak q))\cong \Hom(\frac{R}{\frak a+\frak q},E(R/\frak m))$. Since $\frac{R}{\frak a+\frak q}$ has finite length, $\Hom(\frac{R}{\frak a+\frak q},E(R/\frak m))$ has finite length and since $(0:_{E(R/\frak m)}\frak q)$ is artinian, it follows from [M1, Theorem 1.6] that $(0:_{E(R/\frak m)}\frak q)\in\cS_{\cM.\frak a-cof}$, and hence the assumption on $\frak q$ implies that $E(R/\frak m)\in\cS_{\cM.\frak a-cof}$. Then $\Hom_R(R/\frak a,E(R/\frak m))\cong E_{R/\frak a}(R/\frak m)$ is finitely generated and so $\dim R/\frak a=0$ which contradicts the fact that ${\rm ht \frak p}=\dim R-1$.
In order to prove the second assertion, assume that $\cS_{\cM.\frak a-cof}$ satisfies $C_{\frak a}$ condition. Replacing $\frak p$ by $\frak a$ in the argument mentioned in the first assertion, we deduce that $H_{\frak a}^i(R)\in\cS_{M.\frak a-cof}$ for each $i\geq 0$. Therefore, there exists a $\frak p\in\Min\frak a$ such that ht$\frak a=$ht$\frak p>0$ and $H_{\frak p R_{\frak p}}^{{\rm ht}\frak a}(R_{\frak p})$ is a finitely generated $R_{\frak p}$-module which contradicts the non-vanishing theorem (see [BS, 6.1.7 Exercise]). 
\end{proof}


\section{Local cohomology of weakly Laskerian modules}

In this section, we study the weakly Laskerian modules and their local cohomology modules. We start with a result which characterizes the local cohomology modules from below are in $\cW$.
\medskip
\begin{Theorem}
Let $\frak a$ be an ideal of local ring $R$ of dimension one, let $n\in\bbN$ and let $M$ be an $R$-module. Then $\Ext^i_R(R/\frak a,M)\in\cW$ for all $0\leq i\leq n$ if and only if $H_{\frak a}^i(M)\in\cW$ for all $0\leq i\leq n$.
\end{Theorem}
\begin{proof}
By virtue of \cref{coo}, the subcategory $\cW$ satisfies $C_{\frak a}$ condition. Now, the assertion follows from [AM, Theorem 2.9]
\end{proof}

For an arbitrary ideal $\frak a$ of $R$, we have the following version of the above theorem without the condition $C_{\frak a}$ for $\cW$.
\medskip

\begin{Theorem}\label{s11}
Let $\frak a$ be an ideal of $R$, let $M$ be an $R$-module and let $n$ be a non-negative number. If $H_{\frak a}^i(M)\in\cW$ for all $0\leq i\leq n$, then $\Ext_R^i(R/\frak a, M)\in\cW$ (this also holds for any Serre subcategory $\cS$ of modules). The converse of implication holds if $\Ext_R^1(R/\frak a,K)\in\cW$ for all proper submodules $K$ of $H_{\frak a}^i(M)$ and for all $0\leq i\leq n$.  
 \end{Theorem}
\begin{proof}
We proceed by induction on $n$. The case $n=0$ is trivial as $\Hom_R(R/\frak a,\Gamma_{\frak a}(M))=\Hom_R(R/\frak a, M)$. For $n>0$, we may assume that $\Gamma_{\frak a}(M)=0$ and so there is an exact sequence of modules $0\To M\To E\To Q\To 0$ such that $E$ is injective and $\Gamma_{\frak a}(E)=0$. Now, using the induction hypothesis for $Q$, the result follows straightforwardly. For the converse, if $n=0$ and $K$ is any submodule of $\Gamma_{\frak a}(M)$, the exact sequence $0\To K\To \Gamma_{\frak a}(M)\To\Gamma_{\frak a}(M)/K\To 0$ induces an exact sequence $0\To(0:_K\frak a)\To(0:_{\Gamma_{\frak a}(M)}\frak a)\To(0:_{\Gamma_{\frak a}(M)/K}\frak a)\To \Ext_R^1(R/\frak a,K)$ and hence $(0:_{\Gamma_{\frak a}(M)/K}\frak a)\in\cW$. Since $(0:_{\Gamma_{\frak a}(M)/K}\frak a)$ is an essential submodule of $\Gamma_{\frak a}(M)/K$, the set $\Ass \Gamma_{\frak a}(M)/K=\Ass(0:_{\Gamma_{\frak a}(M)/K}\frak a)$ is finite. Consequently $\Gamma_{\frak a}(M)\in\cW$. Assume, inductively,  that $n>0$ and the result has been prove for all values smaller than $n$. Since $\cW$ is Serre, applying the functor $\Hom_R(R/\frak a,-)$ to he exact sequence $0\To \Gamma_{\frak a}(M)\To M\To M/\Gamma_{\frak a}(M)\To 0$, it is straightforward to see that $\Ext_R^i(R/\frak a, M/\Gamma_{\frak a}(M))\in\cW$ for all $0\leq i\leq n$ and so we may assume that $\Gamma_{\frak a}(M)=0$. Then there exists an exact sequence $0\To M\To E\To Q\To 0$ such that $E$ is injective and $\Gamma_{\frak a}(E)=0$. Thus, there are the isomorphisms $H_{\frak a}^i(Q)\cong H_{\frak a}^{i+1}(M)$ and $\Ext_R^i(R/\frak a,Q)\cong\Ext_R^{i+1}(R/\frak a,M)$ for all $0\leq i\leq n-1$. Furthermore, the assumption implies that $\Ext_R^1(R/\frak a,K)\in\cW$ for all proper submodules $K$ of $H_{\frak a}^i(Q)$ and all $0\leq i\leq n-1$; and hence the induction hypothesis implies that $H_{\frak a}^{i}(M)\cong H_{\frak a}^{i-1}(Q)\in\cW$ for all $1\leq i\leq n$.
\end{proof}

The following theorem about local cohomology from above, is a generalization of [AM, Theorem 3.1] for weakly Laskerian modules.

\medskip

\begin{Theorem}\label{serw}
Let $\frak a$ be an ideal of $R$, let $M$ be a weakly Laskerian $R$-module, let $\cS$ be a Serre subcategory of $R$-modules and let $n$ be a positive integer. Then $H_{\frak a}^i(M)\in\cS$ for all $i\geq n$ if and only if $H_{\frak a}^i(R/\frak p)\in\cS$ for all $\frak p\in\Supp M$ and all $i\geq n$ whenever one of the following conditions holds.

${\rm (i)}$ $n\geq 3$; 

${\rm (ii)}$ If $(R,\frak m)$ is a local ring, $\cS$ satisfies $C_{\frak a}$ condition and $\Hom_R(R/\frak a,M),\Ext_R^1(R/\frak a,M)$ are finitely generated.
\end{Theorem}
\begin{proof}
Since $M$ is weakly laskerian, there exists a finitely generated submodule $N$ of $M$ and a module $F$ of finite support such that $0\To N\To M\To F\To 0$ is exact. We note that $\dim F\leq 1$. To prove the assertion by condition (i), applying the functor $\Gamma_{\frak a}(-)$, there is an isomorphism $H_{\frak a}^i(N)\cong H_{\frak a}^i(M)$ for each $i\geq 3$. If $H_{\frak a}^i(M)\in\cS$ for all $i\geq n$, it follows from [AM, Theorem 3.1] that $H_{\frak a}^i(R/\frak p)\in\cS$ for all $i\geq n$ and all $\frak p\in\Supp N$. If $\frak p\in\Supp F$, then $\dim R/\frak p\leq 1$ and so $H_{\frak a}^i(R/\frak p)=0$ for all $i\geq n$. The converse of implication also follows from [AM, Theorem 3.1]. To do by condition (ii), in view of the above exact sequence and the assumption, $\Hom_R(R/\frak a,F)$ and $\Ext_R^1(R/\frak a,F)$ are finitely generated. If we consider $S=R/\Ann(F)$, then $\dim S/\frak a S \leq 1$. Consider the Grothendieck spectral sequence $$E_2^{p,q}:=\Ext_S^p(\Tor_q^R(S,R/\frak a),F)\Rightarrow H^{p+q}=\Ext_R^{p+q}(R/\frak a,F).$$ We observe that the module $\Hom_S(S/\frak aS,F)\cong \Hom_R(R/\frak a,F)$ is finitely generated. 
Since $\Supp_S\Tor_q^R(S,R/\frak a)\subseteq V(\frak aS)$ and $\Tor_q^R(S,R/\frak a)$ is a finitely generated $S/\frak aS$-module, we conclude that $E_2^{0,q}$ is finitely generated for all $q$. But there are the isomorphisms $E_2^{1,0}\cong E_3^{1,0}\cong\dots\cong E_{\infty}^{1,0}$. Moreover, there are the following filtration 
$$0=\Phi^{2}H^1\subset\Phi^1H^1\subset \Phi^0 H^1\subset H^1$$  such that $E_{\infty}^{1,0}\cong \Phi^1H^1/\Phi^{2}H^1 =\Phi^1H^1$ is a submodule of $H^1=\Ext_R^1(R/\frak a,F)$; and hence it is finitely generated. Now, it follows from [BNS, Theorem 2.5 (i)] that $H_{\frak aS}^i(F)$ is $\frak aS$-cofinite for all $i\geq 0$ and so by virtue of [DM, Proposition 2], $H_{\frak a}^i(F)$ is $\frak a$-cofinite for all $i\geq 0$. Thus $\Hom_R(R/\frak a,H_{\frak a}^1(F))$ is finitely generated and since $\Supp H_{\frak a}^1(F)\subseteq V(\frak m)$, the module $\Hom_R(R/\frak a,H_{\frak a}^1(F))$ has finite length and the fact that $R/\frak m\in\cS$ forces $\Hom_R(R/\frak a,H_{\frak a}^1(F))\in\cS$; and hence $H_{\frak a}^1(F)\in\cS$ as $\cS$ satisfies $C_{\frak a}$ condition. For the case $n=2$, viewing the exact sequence $H_{\frak a}^1(F)\To H_{\frak a}^2(N)\To H_{\frak a}^2(M)\To 0$ and using a similar proof of (i), the assertion follows. For the case $n=1$, assume that $H_{\frak a}^i(M)\in\cS$ for all $i\geq 1$. Since $H_{\frak a}^1(F)\in\cS$, we have $R/\frak m\in\cS$. For every $\frak p\in\Supp F$, if $\frak a\subseteq \frak p$, we have $H_{\frak a}^i(R/\frak p)=0$ for all $i\geq 1$. If $\frak a\nsubseteq\frak p$, there exists $x\in\frak a\setminus \frak p$ and so an exact sequence $0\To R/\frak p\stackrel{x.}\To R/\frak p\To \frac{R}{\frak p+xR}\To 0$ so that $\frac{R}{\frak p+xR}$ has finite length. Since $R/\frak m\in\cS$, we have $\frac{R}{\frak p+xR}\in\cS$. This implies that $(0:_{H_{\frak a}^1(R/\frak p)}x)\in\cS$ so that $(0:_{H_{\frak a}^1(R/\frak p)}\frak a)\in\cS$. Then the assumption on $\cS$ implies that $H_{\frak a}^1(R/\frak p)\in\cS$. We also observe that $H_{\frak a}^i(R/\frak p)=0$ for all $i\geq 2$. Now, we want to show that $H_{\frak a}^{i}(R/\frak p)\in\cS$ for each $\frak p\in\Supp M$ and each $i\geq 1$. By a descending induction, we may prove that $H_{\frak a}^{1}(R/\frak p)\in\cS$ for each $\frak p\in\Supp M$. Assume that $\frak p\in\Supp M$ is maximal of those $\frak p\in\Supp M$ such that $H_{\frak a}^{1}(R/\frak p)\notin\cS$. The previous argument implies that $\frak p\notin\Supp F$ so that $M_{\frak p}=N_{\frak p}$ is finitely generated. Then Nakayama lemma implies that $M_{\frak p}/\frak p M_{\frak p}\neq 0$ and a similar proof of [Bo, Chap, (ii), 4, Proposition 20] implies that $\Hom_R(M,R/\frak p)\neq 0$. Applying a similar argument mentioned in the proof of [AM, Theorem 3.1], we deduce that $H_{\frak p}^{1}(R/\frak p)\in\cS$ which is a contradiction. Conversely assume that $H_{\frak a}^i(R/\frak p)\in\cS$ for all $i\geq 1$ and all $\frak p\in\Supp M$. It follows from [AM, Theorem 3.1] that $H_{\frak a}^i(N)\in\cS$ for all $i\geq 1$.  Now, the exact sequence $$H_{\frak a}^1(N)\To H_{\frak a}^1(M)\To H_{\frak a}^1(F)\To H_{\frak a}^2(N)\To H_{\frak a}^2(M)\To 0$$ and the previous argument imply that $H_{\frak a}^1(M),H_{\frak a}^2(M)\in\cS$. Moreover, the isomorphisms $H_{\frak a}^i(N)\cong H_{\frak a}^i(M)$ for all $i\geq 3$ imply that $H_{\frak a}^i(M)\in\cS$ for all $i\geq 3$. 
\end{proof}

Given a Serre subcategory $\cS$ of $R$-modules, an ideal $\frak a$ of $R$ and a positive integer $n$, we denote by $\cS_n(\frak a)$ the subcategory of all modules $M$ such that $H_{\frak a}^i(M)\in\cS$ for all $i\geq n$. Clearly $\cS_1(\frak a)\subseteq\cS_2(\frak a)\subseteq \dots$ and $\cS_i=R$-Mod for all $i>\dim R$, where $R$-Mod is the category of $R$-modules. It follows from [AM, Theorem 2.8] that $\cS_n(\frak a)\cap\cN$ is a Serre subcategory for each $n\geq 1$. Furthermore, \cref{serw} implies that $\cS_n(\frak a)\cap\cW$ is a Serre subcategory for each $n\geq 3$.

\medskip

\begin{Proposition}
Let $\frak a$ be an ideal of $R$. Then the following conditions hold.
 
${\rm (i)}$ There is a chain $\bbM[\cS_1(\frak a)\cap\cN]\supseteq\bbM[\cS_2(\frak a)\cap\cN]\supseteq\dots$.

${\rm (ii)}$ If $(R,\frak m)$ is a local ring, then there is a chain $\bbM[\cS_3(\frak a)\cap\cW]\supseteq\bbM[\cS_4(\frak a)\cap\cW]\supseteq\dots$.
\end{Proposition}      
\begin{proof}
(i) It suffices to show that $\bbM[\cS_1(\frak a)\cap\cN]\supseteq\bbM[\cS_2(\frak a)\cap\cN]$. Assume that $\frak p\in\bbM[\cS_2(\frak a)\cap\cN]$ and $M$ is an $R$-module such that $\Gamma_{\frak p}(M)=M$ and $(0:_M\frak p)\in\cS_1(\frak a)\cap\cN$. Then$(0:_M\frak p)\in\cS_2(\frak a)\cap\cN$ so that $M\in\cS_2(\frak a)\cap\cN$. It follows from [AM, Theorem 3.1] that $H_{\frak a}^i(R/\frak q)\in\cS$ for all $\frak q\in\Supp(0:_M\frak p)=\Supp M$ and all $i\geq 3$. Using again  [AM, Theorem 3.1], we deduce that $H_{\frak a}^i(M)\in\cS$ for each $i\geq 1$. (ii) It suffices to show that $\bbM[\cS_3(\frak a)\cap\cW]\supseteq\bbM[\cS_4(\frak a)\cap\cW]$. If $\frak p\in\bbM[\cS_4(\frak a)\cap\cW]$ and $M$ is an $R$-module such that $\Gamma_{\frak p}(M)=M$ and $(0:_M\frak p)\in\cS_3(\frak a)\cap\cW$. Then we conclude that $M\in\cS_4(\frak a)\cap\cW$ and so there is an exact sequence $0\To N\To M\To F\To 0$ such that $N$ is finitely generated and $F$ has finite support. Since $\Supp F\subseteq V(\frak p)$, we have $\Supp F=\{\frak p,\frak m\}$; and hence $\Supp N=\Supp M$. On the other hand, $H_{\frak a}^i((0:_M\frak p))\in\cS$ for each $i\geq 3$ and so it follows from \cref{serw} that $H_{\frak a}^i(R/\frak q)\in\cS$ for all $\frak q\in\Supp (0:_M\frak p)=\Supp N$ for each $i\geq 3$. Thus [AM, Theorem 2.8] implies that $H_{\frak a}^i(N)\in\cS$ for each $i\geq 3$; and consequently the exact sequence $H_{\frak a}^i(N)\To H_{\frak a}^i(M)\To 0$ implies that $H_{\frak a}^i(M)\in\cS$ for each $i\geq 3$.   
\end{proof}

The rest of this section is devoted to the cofiniteness of local cohomology modules of weakly Laskerian modules.  For an ideal $\frak a$ of $R$, an $R$-module $M$ is said to be $\frak a$-{\it cofinite} if $\Supp M\subseteq V(\frak a)$ and
$\Ext_R^i(R/\frak a, M)$ is finitely generated for all integers $i\geq 0$. For the first local cohomology, we have the following proposition.

\medskip
\begin{Proposition}\label{pp}
If $\frak a$ is an ideal of $R$ and $M$ is a weakly Laskerian module and $\Ext_R^i(R/\frak a,M)$ is finitely generated for $i=0,1$, then $\Gamma_{\frak a}(M)$ is $\frak a$-cofinite. 
\end{Proposition}
\begin{proof}
The assumption implies that $\Ext_R^i(R/\frak a,\Gamma_{\frak a}(M))$ is finitely generated for $i=0,1$. Since $\Gamma_{\frak a}(M)$ is weakly Laskerian, there is an exact sequence $0\To N\To\Gamma_{\frak a}(M)\To F\To 0$ such that $N$ is finitely generated and $\dim F\leq 1$. The previous argument implies that $\Ext_R^i(R/\frak a,F)$ is finitely generated for $i=0,1$ and $\Supp F\subseteq V(\frak a)$. Therefore, it follows from [BNS, Lemma 2.2] that $F$ is $\frak a$-cofinite and consequently $\Gamma_{\frak a}(M)$ is $\frak a$-cofinite.  
\end{proof}

\medskip
For the $i$-th local cohomology of weakly Laskerian modules, we have the following theorem. 
\begin{Theorem}\label{ll}
Let $\frak a$ be an ideal of $R$ such that $\dim R/\frak a=1$ and let $M$ be a weakly laskerian $R$-module such that $\Ext_R^i(R/\frak a,M)$ is finitely generated for $i=0,1$. Then $H_{\frak a}^i(M)$ is $\frak a$-cofinite for all $i\geq 0$. 
\end{Theorem}
\begin{proof}
There is an exact sequence $0\To N\To M\To F\To 0$ such that $N$ is finitely generated and $F$ has finite support; and hence $\dim F\leq 1$. The assumption implies that $\Ext_R^i(R/\frak a,F)$ is finitely generated for $i=0,1$ and hence it follows from [BNS, Theorem 2.5(i)] that $\Ext_R^i(R/\frak a,F)$ is finitely generated for all $i\geq 0$. It now follows from [BNS, Theorem 2.5(ii)] that $H_{\frak a}^i(M)$ is $\frak a$-cofinite for all $i\geq 0$. 
\end{proof}

\medskip
The following theorem  generalizes [NS, Theorem 3.7] for the weakly Laskerian modules.

\begin{Theorem}\label{dimtwo}
Let $R$ be a local ring, let $\frak a$ be an ideal of $R$ such that $\dim R/\frak a=2$ and let $M$ be a weakly Laskerian $R$-module such that $\Ext_R^i(R/\frak a,M)$ is finitely generated for $i=0,1$. Let $n$ be a non-negative integer. Then the following conditions are equivalent:

${\rm (i)}$ $\Hom_R(R/\frak a, H_{\frak a}^i(M))$ is finitely generated for all $i\leq n$.

${\rm (ii)}$ $H_{\frak a}^i(M)$ is $\frak a$-cofinite for all $i<n$.   
\end{Theorem}
\begin{proof}
If $n=0$, there is nothing to prove and so we may assume that $n>0$. Since $M$ is weakly Laskerain, there exists an exact sequence of modules $0\To N\To M\To F\To 0$ such that $N$ is finitely generated and $\dim F\leq 1$. The assumption implies that $\Ext_R^i(R/\frak a,F)$ is finitely generated for $i=0,1$; and hence it follows from [BNS, Theorem 3.4] that $\Ext_R^i(R/\frak a,F)$ is finitely generated for all $i\geq 0$. Therefore $\Ext_R^i(R/\frak a,M)$ is finitely generated for all $i\geq 0$. Now the result follows by using [NS, Theorem 3.7].
\end{proof}
\medskip
\begin{Corollary}
Let $R, \frak a,M$ be as in \cref{dimtwo}. Then $\Hom_R(R/\frak a,H_{\frak a}^1(M))$ is finitely generated.
\end{Corollary}
\begin{proof}
It is straightforward by \cref{pp} and \cref{dimtwo}.
\end{proof}
\medskip

\begin{Corollary}\label{dim3}
If $R$ is a local ring of dimension $\leq 3$, then \cref{dimtwo} holds for any ideal $\frak a$ of $R$. 
\end{Corollary}
\begin{proof}
Assume that $\frak a$ is an ideal of $R$ and $M$ is a weakly Laskerian module such that $\Ext_R^i(R/\frak a,M)$ is finitely generated for $i=0,1$. By virtue of \cref{pp}, the module $\Gamma_{\frak a}(M)$ is $\frak a$-cofinite; and hence we may assume that $\Gamma_{\frak a}(M)=0$. On the other hand, there exists a positive integer $t$ such that $\Gamma_{\frak a}(R)=(0:_R\frak a^t)$. If $\frak a^t=0$, there is nothing to prove; otherwise take $S=R/\Gamma_{\frak a}(R)$. Clearly $M$ is an $S$-module and $\dim S/\frak aS\leq 2$. By the a similar proof mentioned in \cref{serw}, we deduce that $\Hom_S(S/\frak aS,M)=0$ and $\Ext^1_S(S/\frak aS,M)$ is finitely generated. For each $i$, there is an isomorphism  $H_{\frak a}^i(M)\cong H_{\frak aS}^i(M)$ and it follows from [DM, Proposition 2] that $H_{\frak a}^i(M)$ is $\frak a$-cofinite if and only if $H_{\frak aS}^i(M)$ is $\frak aS$-cofinite and $\Hom_R(R/\frak a, H_{\frak a}^i(M))\cong\Hom_S(S/\frak aS,H_{\frak a}^i(M))$. Furthermore, it is clear that $M$ is weakly Laskerian $S$-module. If $\dim S/\frak aS=1$, the result follows by using \cref{ll}. If $\dim S/\frak aS=2$, the result follows by \cref{dimtwo}.  
\end{proof}



\end{document}